\documentclass[10pt, a4paper]{amsart}
\input xy   
\xyoption{all}
\numberwithin{equation}{section}
\swapnumbers

\theoremstyle{plain}
  \begingroup
        \newtheorem{theorem}[equation]{Theorem}
        
        \newtheorem{proposition}[equation]{Proposition}
        \newtheorem{corollary}[equation]{Corollary}
        
        \newtheorem{facts}[equation]{Facts}
        \newtheorem{assumption}[equation]{Assumption}
        \newtheorem{remark}[equation]{Remark}
	\newtheorem{definition}[equation]{Definition}
        
        \newtheorem{sinnadaitalica}[equation]{}
  \endgroup

\theoremstyle{definition}
  \begingroup

        \newtheorem{sinnadastandard}[equation]{}
        \newtheorem{observation}[equation]{Observation}
  \endgroup

\newcommand{\mr}[1]{\buildrel {#1} \over \longrightarrow}
\newcommand{\ml}[1]{\buildrel {#1} \over \longleftarrow}

\newcommand{\rimply}{\Rightarrow}

\newcommand{\cc}{\mathcal}
\newcommand{\bb}{\mathbb}
\newcommand{\mbf}{\mathbf}

\newcommand{\siff}{\Leftrightarrow}
\newcommand{\mono}{\hookrightarrow}
\newcommand{\mmr}[1]{\buildrel {#1} \over \hookrightarrow}

\newcommand{\pa}{\cdot}

\begin{document}

\title{The fundamental progroupoid of a general topos}

\author{Eduardo J. Dubuc}

\begin{abstract}
It is well known  that the category of covering projections (that is,
  locally constant objects) of a locally connected topos is equivalent
  to the classifying topos of a strict progroupoid
  (or, equivalently, 
  a localic prodiscrete groupoid), the 
  \emph{fundamental progroupoid}, and that this
  progroupoid represents first degree cohomology. In this paper we
  generalize these results to an arbitrary topos. The fundamental
  progroupoid is now a localic progroupoid, and can not be replaced by
  a localic groupoid. The classifying topos in not any more a Galois
  topos. 
Not all locally constant objects can be considered as covering
projections. The
key  contribution of this paper is a novel definition of
covering projection for a general topos, which coincides with the
usual definition when the topos is locally connected. The results in
this paper were presented in a talk at the Category Theory Conference,
Vancouver July 2004.
\end{abstract}

\maketitle

{\sc introduction.} 
It is well known that if $\cc{E}$ is a locally connected topos then
the category of covering projections (that is, locally constant
objects) is equivalent to the classifying topos of a strict progroupoid (or,
  equivalently, a localic prodiscrete groupoid), the
  \emph{fundamental progroupoid} $\pi(\cc{E})$, and that this
  progroupoid represents first degree cohomology. In this paper we
  generalize these results to an arbitrary topos. 

The subject that concern us here was developed (in the context of
grothendieck topos over $\cc{S}et$) for a pointed locally connected topos by
Grothendieck-Verdier in a series of commented exercises in Expose
IV of the SGA4 \cite{G2}, and by Artin-Mazur \cite{AM}. Later Moerdijk
\cite{M2} treats the subject over a general base topos $\cc{S}$, and
replace \emph{progroups} by \emph{prodiscrete localic groups}.
Bunge \cite{B} does the unpointed case and works with \emph{prodiscrete
localic groupoids}. See also Bunge-Moerdijk \cite{BM}, and the Appendix in
\mbox{Dubuc \cite{D2}} for a 
resume of this theory.

The salient feature of the theory is that covering projections are
considered as a \emph{full} subcategory of the topos, and this fact is
essential in the proofs of the validity of the statements. Covering
projections can not be considered as a full subcategory when the topos
is not locally connected.

The principal source of inspiration for our work was the paper of
\mbox{Hernandez-Paricio 
\cite{LJHP},} where he treats successfully the case of non locally connected
topological spaces. There it is possible to see that there is a descent
datum underneath the notion of covering projection, and that this datum
has to be taken into account in the definition of covering
projection. We can see an
implicit situation of classical 
topological descent as described in the 
introduction to ``Categories Fibrees et Descente'', \cite{G1}, Expose VI.
Once the descent datum is made explicit, the category of covering
projections of a topological space trivialized by a (fix)
covering is, by its very definition, the 
classifying topos of a discrete groupoid, and  
this groupoid can be explicitly constructed as 
the free category over the nerve of the covering. The assignemment of
this groupoid is functorial on the filtered poset of covering sieves
and determines the fundamental progroupoid of the space. We explain all this
in \mbox{section 1.}    

Given any topos, there is no problem to
construct the topos of locally
constant objects trivialized by a (fix) cover. The problem
is that when the topos is not locally connected, the resulting topos
is not atomic   
because it fails to be both locally connected and boolean. The
situation here is not as simple as in the case of a topological space,
 but nevertheless the work in \cite{LJHP} gave us the clue to
 unravel  this situation. Namely, \emph{even though Hernandez-Paricio
 covering projections are not locally connected topological spaces,
 their set of connected components can be constructed}.

In section 2 we prepare the ground for our work by explicitly
establishing an equivalence between the the usual definition of locally
constant object and a certain descent datum. Section 3 contains the
key  contribution of this paper, which is a novel definition of
covering projection for a general topos (definition
\ref{coveringprojectionobject}). When the topos is locally connected,
every locally constant object is a covering projection in our
sense. We construct the topos of 
covering projections and show that it is atomic (theorem
\ref{atomic}).
The resulting localic groupoid (theorem \ref{U_fundamental}) 
appears to be the first genuine application of Joyal-Tierney results 
(\cite{JT} VIII, 3. \mbox{Theorem 1)} to the galois theory of locally constant
objects. Here for the first time a non
prodiscrete localic groupoid appears in this theory, as well as an
atomic topos which is not a Galois topos. In sections 4, 5 and 6 we
show that our notion of covering projection is well behaved and
adequate to a treatment with inverse limit of topoi techniques.
We construct the category of all covering projections, the topos it
generates, and the fundamental (in this case localic)
progroupoid. We show an equivalence between the classifying topos of
this progroupoid and the topos of covering projections
\mbox{(theorem \ref{fundamentaltopos}).} Finally, in section 7 we prove that
torsors (for a discrete group) are covering projections in our sense,
and that the fundamental localic progroupoid represents first degree
cohomology (theorem \ref{representstorsors}).  
         
\vspace{1ex}

\emph{Comparison between the locally connected, spacial and general
  cases.}

In the case of locally connected topoi, given a (fix) cover, the
points of the 
topos of covering projections are essential, and the
corresponding groupoid is an ordinary (discrete) groupoid.  This
determines a fundamental ordinary progroupoid.  The transition morphisms are
surjective on triangles, fact that allows to replace this progroupoid
by a prodiscrete localic groupoid. Equivalent to this, the transition
morphisms between the topoi are connected, and this implies that the
classifying topos is a Galois Topos. 

In the case of a non locally connected topological space, given a
(fix) cover, the corresponding 
groupoid is still discrete, and we still have an ordinary fundamental
progroupoid. However, it can not be replaced by a localic groupoid
because the transition morphisms are not surjective on triangles, or ,
equivalently, the transition morphisms between the topoi are not
connected. The classifying topos is not any more a Galois topos.

In the case of a general topos, given a (fix) cover, the points of the
topos of covering projections are not essential, and the corresponding
groupoid is a localic (non discrete) groupoid. This determines a
fundamental localic progroupoid. Neither the topoi in the system nor
the classifying topos are Galois.

\vspace{1ex} 

{\sc context.} 
Throughout this paper $\cc{S} = Sets$ denotes the topos of sets. However, we
argue in a way that should be valid if $\cc{S}$ is
an arbitrary grothendieck topos, but let the interested reader to verify
this. All topoi $\cc{E}$ are assumed to be grothendieck topoi
(over $\cc{S}$), the structure map will be denoted by 
$\gamma: \cc{E} \to \cc{S}$ in all cases. 

Recall that a geometric morphism  $\cc{E} \mr{f} \cc{F}$ is said to be 
\emph{essential} if the inverse image functor $f^{*}$ has
itself a left adjoint $f_{\,!}$, \emph{locally connected} when
$f_{\,!}$ is $\cc{F}$-indexed,   
 \emph{connected} if
the inverse image functor $f^{*}$ is full and faithful, and
\emph{atomic} if $f^{*}$ is logical. A topos is said to be
locally connected, connected, or atomic, when the structure morphism
is so. A topos is atomic if and only if it is locally connected and
boolean. We refer to
\cite{G2} Expose IV, 4.3.5, 4.7.4, 7.6 and 8.7., \cite{M1},
and \cite{JT}.

\tableofcontents

\section{Covering projections of an arbitrary topological space}
\label{cpts} 

Topologists have dealt successfully
with covering projections of non locally connected topological spaces.
 In their work, the descent
data underneath the notion of covering projection has to be made
explicit in one way or another (see 
\cite{LJHP} and references therein).
This shows that we are in face of a situation of classical
topological descent as described in the 
introduction to ``Categories Fibrees et Descente'', \cite{G1}, Expose VI.

Given a topological space $B$, when $B$ is not locally
connected, covering projections can not
be  considered as local homeomorphisms of a particular type, but
should be considered as local homeomorphisms
with an added structure. This is reflected by the fact that not all
continuous maps between the underlying sheaves are admitted, but only
those that preserve the trivialization structure. This determines a
category (topos) $\cc{P}_\cc{U}$, a geometric morphism
 $Sh(X) \to \cc{P}_\cc{U}$, with faithful but not full inverse
image $\cc{P}_\cc{U} \to Sh{X}$) (where $Sh(B)$ is the topos of
  sheaves over $B$), and a surjective point
$\cc{S}_{/I} \to  \cc{P}_\cc{U}$, which is not (contrary to the case of
  a locally connected space)  of effective descent. 
 
Not all
locally constant sheaves over $B$ should be admited as covering
projections. Consider a sheaf $X \to B$ split by a
cover $\cc{U} = \{U_i\}_{i \in I}$, $U_i \subset B$, by means of homeomorphisms
 \mbox{$S_{i} \times  U_{i} \mr{\theta_{i}} X|_{U_i}$}. 
Given $i,\, j,\; U_i \cap U_j \neq \emptyset$, there is an
induced homeomorphism 
\mbox{$S_i \times (U_i \cap U_j)
\mr{\theta_{i}} X|_{U_i \cap U_j}  
\mr{\theta_{j}^{-1}} S_j \times (U_i \cap U_j)$} over $ U_i \cap
U_j$. The following definition is essentially definition 2.1 in \cite{LJHP}.
\begin{definition} \label{coveringprojection}
A covering projection split by $\cc{U}$ is a locally constant sheaf $X$ such that the bijections
between the fibers $S_i  \times \{x\}  \to  S_j \times
\{x\}$ are given by the same function for all the points 
$x \in U_i \cap U_j$. That is, $\forall x, y \in U_i \cap U_j$, 
\mbox{$\theta_{j}^{-1} \circ \theta_{i}(-, \, x) = \theta_{j}^{-1} \circ
\theta_{i}(-, \, y)$.}
\end{definition}
%
%
Covering projections $X \to B$ trivialized by $\cc{U}$
define a full
subcategory $\cc{D}_\cc{U} \subset \cc{P}_\cc{U}$,  $\cc{D}_\cc{U}$ is
a topos, and there is a (connected) geometric morphism
$\cc{P}_\cc{U} 
\to \cc{D}_\cc{U}$ with inverse image given by the full inclusion. Thus
$\cc{D}_\cc{U}$ is equipped with a surjective point 
\mbox{$\cc{S}_{/I} \to \cc{P}_\cc{U} \to \cc{D}_\cc{U}$} 
which is of effective descent. 
In fact, the collection of
homeomorphisms \mbox{$\theta_{j}^{-1} \circ \theta_{i}: S_i \times
(U_i \cap U_j) \to S_j \times (U_j \cap U_i)$} \mbox{defines} a
situation of classical (topological) descent \cite{G1} Expose VI.
The condition in definition \ref{coveringprojection} means that there are bijections 
$\lambda_{j\,i}: S_i \to  S_j$ which induce the composite
homeomorphisms $\theta_{j}^{-1} \circ \theta_{i}$. The sheaf
$X$ together with the trivialization $\{\theta_i\}_{i \in I}$ can be recovered  by descent from the family of
topological spaces $X_i = S_i \times U_i$ and the bijections
$\lambda_{j\,i}: S_i \to  S_j$. 
The cover $\cc{U}$
determines a simplicial set $\cc{U}_\bullet$ (the Cech nerve) whose
$n$-simplexes are given by   
$\cc{U}_n = \{(i_0, i_1, \ldots i_n) \;|\;  U_{i_0} \cap  U_{i_1}
  \ldots 
\cap U_{i_n} \;\neq\; \emptyset\}$ (notice that $\cc{U}_0 = I$). In
turn, this determines a simplicial topos $\cc{S}_{/\cc{U}_\bullet}$ by slicing.
The family of
bijections $\lambda_{j\,i}$ is exactly a
\emph{descent datum} on the object 
$S \to I$ in the topos $\cc{S}_{/I}$. The topos $\cc{D}_\cc{U}$ is
(equivalent to) the descent topos 
$\cc{S}_{\cc{U}_\bullet} \mr{}\cc{D}_{\cc{U}}$, and as such, the
morphism $\cc{S}_{/I} \to \cc{D}_\cc{U}$ is  of effective descent.

\emph{We have then that, once the descent datum underneath the notion
  of covering projection is made explicit, the category of covering
projections ($B$ locally connected or not) trivialized by a (fix)
covering $\cc{U}$ is, by its very definition, the 
classifying topos of a discrete groupoid $\pi\cc{U}$ (whose objects are the index
set of the cover).} 
This groupoid can be explicitly constructed as 
  the free category over the nerve of the covering.   

\vspace{1ex}

The collection of faithful functors $\cc{D}_\cc{U} \to Sh(B)$ form a
  cone over the category of refinements, and this allows the
  construction of the category of all
  covering projections $Cp(B) \to Sh(B)$ as the colimit of the system
  of categories  
$\cc{D}_\cc{U}$ \cite{LJHP}. The system of groupoids $\pi \cc{U}$, with $\cc{U}$
  running over the filtered poset of covering sieves, determines a
  progroupoid, whose category of actions (as defined in
  \cite{LJHP}) is equivalent to $Cp(B)$.

\vspace{1ex}

When the space $B$ is locally connected, the condition on definition
\ref{coveringprojection} is vacuous, $\cc{P}_\cc{U} = \cc{D}_\cc{U}$,
and this theory yields the classical galois theory of locally
connected topological spaces.

\vspace{2ex}

\section{Locally constant objects and descent data}
\label{lco} 

By a \emph{cover} $\cc{U} = \{U_{i}\}_{i \in I}$ in a topos 
$\cc{E} \mr{\gamma} \cc{S}$ we mean an epimorphic family
$U_{i} \to 1$ in $\cc{E}$, $ I \in \cc{S}$. As usual, this is
alternative notation for a map 
$\zeta \colon U \to {\gamma}^{*}I$, with $U \to 1$ epimorphic. Notice
that covers are $3$-tuples $\cc{U} = (U, I, \zeta)$. 
\begin{assumption} \label{nonempty}
We assume that $U_i \neq \emptyset \; \forall i \in I$.
\end{assumption}   
   
The concept of \emph{locally constant object} is a direct translation
into the topos context of the classical notion of 
\emph{covering projection} (\cite{S} Ch. 2, Sec. 1). It is defined in
SGA4 Expose IX, 2.0. 
%
%
\begin{definition} \label{trivialization} Given a topos $\cc{E} \mr{\gamma}
  \cc{S}$ and a cover $\cc{U} = \{U_{i}\}_{i \in I}$, a
  \emph{trivialization} of an object $X \in \cc{E}$ is a family of
  sets $\{S_{i}\}_{i \in I}$ together with isomorphisms in 
$\cc{E}$,
\mbox{$\{\gamma^{*}S_{i} \times U_{i} \mr{\theta_{i}}  X \times
U_{i}\}_{i \in I}$} over $U_i$. Alternatively, it is an arrow
$S \to I$ in $\cc{S}$, and an isomorphism $\theta \colon {\gamma}^{*}S
\times_{{\gamma}^{*}I} U \to X \times U$ over $U$. We say that $X$ is
  $\cc{U}$-split by the trivialization.
\end{definition}

Trivializations will be denoted by the letter $\theta$ in all cases.

\begin{definition} \label{lc}
An object $X$ in a topos $\cc{E}$ is  \emph{locally
constant} if it is $\cc{U}$-split for some cover $\cc{U}$ in $\cc{E}$.
\end{definition}

In \cite{B} M. Bunge introduces a push-out of topoi
whose underline category is the category of locally
constant objects split by a (fixed) cover, and whose arrows are maps which preserve
the trivializations. 

\begin{definition} [M.Bunge] \label{pushout}
Given a cover 
 $\zeta \colon U \to {\gamma}^{*}I$ in a topos $\cc{E}$, the category
 $\cc{P}_\cc{U}$ of locally constant objects split by $\cc{U}$ is
 given by the following push-out topos:
$$ 
\xymatrix
        {
         \cc{E}_{/U} \ar[r]^{\varphi_\cc{U}} \ar[d]^{\rho_\cc{U}} 
        & \cc{E} \ar[d]^{\upsilon_{\cc{U}}}
        \\
         \cc{S}_{/I} \ar[r]^{p_\cc{U}} 
        & \cc{P}_\cc{U}
        }
$$
where $\rho_{\cc{U}}$ and $\varphi_\cc{U}$ are given by 
$\rho_{\cc{U}}^{*}(S \to I) \,=\, \gamma^{*}S \,
\times_{\gamma^{*}I} \, U$  and
\mbox{$\varphi_{\cc{U}}^{*}(X) \,=\, X \,\times \, U$}. 
By the constructions of push-outs of topoi, the category
$\cc{P}_\cc{U}$ is the following:

\emph{Objects}: $(X, S \to I, \theta)$, is a $3$-tuple  
, with $X \in \cc{E}$, 
$S \to I \,\in\, \cc{S}_{/U}$ and \mbox{$X \times U \mr{\theta} \gamma^{*}S
\times_{\gamma^{*}I} U$}   
an isomorphism over $U$.

\emph{Arrows}: 
$(X, S \to I, \theta) \;\to\;  (Y, T \to I, \theta)$, 
is a pair of morphisms $X \mr{f} Y  \,\in \,\cc{E}$, 
$S \mr{\vartheta} T \,\in\, \cc{S}$ over $I$, $\{S_i \mr{\vartheta_i} T_i\}$,
 compatible with the trivialization data:
$$
\xymatrix
         {
          \gamma^{*}S_i \times U_i 
	         \ar[r]^{\theta_i}
                 \ar[d]^{\gamma^{*}\vartheta_i \times U_i}
          &
          X \times U_i
                 \ar[d]^{f \times U_i}
          \\ 
          \gamma^{*}T_i \times U_i  
                 \ar[r]^{\theta_i}
          &
          Y  \times U_i
         }
$$
The functors $\upsilon_{\cc{U}}^{*} : \cc{P}_\cc{U} \to
\cc{E}$ and $p_{\cc{U}}^{*} : \cc{P}_\cc{U} \to \cc{S}_{/I}$ 
are the projections.
\end{definition} 
It is important to remark that the map $X \mr{f} Y$ completely
determines the function $S \mr{\vartheta} T$. That is, the latter, if it
exists, is unique.  Arrows in  $\cc{P}_\cc{U}$ can be considered as
maps in $\cc{E}$ satisfying a condition. However, on spite of this
uniqueness, we shall also say
that $\vartheta$ lift $f$ into an arrow in $\cc{P}_\cc{U}$. 

\vspace{1ex}

When the
topos $\cc{E}$ is
locally connected, we can assume that $\rho_{\cc{U}}$ is connected and
locally connected (that is, all the objects $U_i$ connected). It follows then
 that the point $\cc{S}_{/I} \mr{p_\cc{U}}
\cc{P}_\cc{U}$ is locally connected and surjective, in particular, of
effective descent \cite{B}. Thus, $\cc{P}_\cc{U}$ is equivalent to the
classifying topos of the (discrete because the point, being locally
connected is representable) groupoid of automorphisms
of $p_\cc{U}$ \cite{JT}, in particular, $\cc{P}_\cc{U}$ is an atomic
topos. It follows also that $\upsilon_{\cc{U}}$ is connected, thus
$\cc{P}_\cc{U} \subset 
\cc{E}$ is a full subcategory via the inverse image functor
$\upsilon_{\cc{U}}^{*}$  \cite{BM}. The
reader can also check the appendix in \cite{D2}, where a all this is
proved ``by hand'', without recourse to the results in \cite{JT}.  
  
%

The known theory for the locally connected case stops here for non locally
connected topoi. It is not any more possible to assume that
$\rho_{\cc{U}}$ is connected and locally connected. Now
$\cc{P}_\cc{U}$ is not any more a full subcategory of $\cc{E}$, and
furthermore, the point
$\cc{S}_{/I} \mr{p_\cc{U}} \cc{P}_\cc{U}$ fails to be of effective
descent.  
We analyzed the situation an found that although $p_\cc{U}$ is
surjective, the problem is that
$\cc{P}_\cc{U}$ is not atomic because it fails to be both locally
connected and boolean. 

\vspace{2ex}

{\bf Trivializations versus descent data}


 Consider a topos $\cc{E} \mr{\gamma} \cc{S}$ and a cover $\cc{U} = (U, I, \zeta)$,
$\zeta \colon U \to {\gamma}^{*}I$, $I \in \cc{S}$, \mbox{$U \in
\cc{E}$}.

\vspace{1ex}

Let $U_\bullet$ (resp. $I_\bullet$) be the simplicial object
(resp. set) whose $n$-simplexes are given by
\mbox{$U_n = U \times U  \times \cdots  U$}
\mbox{(resp. $I_n = I \times I  \times \cdots  I$).}

\vspace{1ex}

Let $\cc{U}_\bullet
\subset I_\bullet$ be the Cech nerve of $\cc{U}$, that is, the
simplicial set 
\mbox{$\cc{U}_n = \{(i_0, i_1, \ldots i_n) \;|\;  U_{i_0} \times  U_{i_1} \ldots
\times U_{i_n} \;\neq\; \emptyset\}$} (notice that $\cc{U}_0 = I$).

\vspace{1ex}

Since the cover will remain fixed in this section, to simplify
notation we shall omit a subindex $\cc{U}$ on the arrows. The map $\zeta :  U \to {\gamma}^{*}I$ induces a morphism
of simplicial objects \mbox{$\zeta_\bullet :  U_\bullet  \to {\gamma}^{*} I_\bullet$}
which factors  through $U_\bullet  \to {\gamma}^{*}
\cc{U}_\bullet \subset {\gamma}^{*} I_\bullet$. Actually, ${\gamma}^{*}
\cc{U}_\bullet$ is the image of $\zeta_\bullet$. We abuse notation and
write $\zeta_\bullet :  U_\bullet  \to {\gamma}^{*}
\cc{U}_\bullet $. This morphism determines a morphism of simplicial topoi
$\rho_{\bullet} : \cc{E}_{/U_\bullet}  \mr{}  \cc{S}_{/
  \cc{U}_\bullet}$. 

\vspace{1ex}

Consider the
truncated simplicial topoi which determine the descent situations (see
for example 3.2 \emph{Descent.} in \cite{M3}):  
\begin{equation} \label{descentsituation}
 \xymatrix@R=7ex
   {
     \cc{E}_{/  U \times U \times U} \ar@<2ex>[r]  \ar[r]
                                \ar@<-2ex>[r] \ar[d]^{\rho_{2}} 
   & \cc{E}_{/  U \times U}  \ar@<2ex>[r]  \ar@<-2ex>[r]
                                              \ar[d]^{\rho_{1}}
   & \cc{E}_{/ U} \ar[l] \ar[r]^{\varphi} \ar[d]^{\rho}
   & \cc{E} \ar@/^/[dr]^{\mu} \ar[d]^{\upsilon}
   \\ 
     \cc{S}_{/ \cc{U}_2} \ar@<2ex>[r]  \ar[r] \ar@<-2ex>[r]
   & \cc{S}_{/ \cc{U}_1 }  \ar@<2ex>[r]  \ar@<-2ex>[r] 
   & \cc{S}_{/I} \ar[l] \ar[r]^{\varrho} \ar`d[r]`[rr]^{\delta}[rr]
   & \cc{P}_\cc{U} \ar[r] 
   & \cc{D}_\cc{U}
  } 
\end{equation}

Where $\cc{S}_{/I} \mr{\delta} \cc{D}_\cc{U}$ is defined to be the
descent topos. It is well known
that the morphism $\cc{E}_{/U} \mr{\varphi} \cc{E}$ is of
effective descent. It follows the  existence of the arrow $\cc{E}
\mr{\mu} \cc{D}_\cc{U}$. $\cc{P}_\cc{U}$ is the push-out topos
\ref{pushout}, and the morphism $\cc{P}_\cc{U} \to \cc{D}_\cc{U}$
follows by the universal property of this push-out.
We shall examine now in more detail diagram \ref{descentsituation}.

\begin{facts} \label{facts}
\end{facts}
\begin{enumerate}
\item \label{phi}
 The morphism $\varphi$ is given by
$$ 
\varphi^{*}(X) \,=\, X \times  U\, :  \;\; \;\; 
\varphi^{*}(X) \,=\, \{X \times U_i\}_{i \in I}.
$$

\vspace{1ex}\item \label{rho}
 The morphism $\rho$ is given by 
$$
\rho^{*}(S \to I)
 \,=\, \gamma^{*}S \, \times_{\gamma^{*}I} \, U \, :\;\; \;\;
\rho^{*}(\{S_i\}_{i \in I}) = \{\gamma^{*}S_i \times U_i\}_{i
 \in I}.
$$
$$
\rho_1^{*}(S \to \cc{U}_1)
 \,=\, \gamma^{*}S \, \times_{\gamma^{*}\cc{U}_1} \, (U \times U) \,:
 \hspace{25ex}
$$
$$ 
\hspace{10ex} \rho^{*}(\{S_{(i,\, j)} \}_{(i,\,j) \in \cc{U}_1}) =
\{\gamma^{*}S_{(i,\,j)} \times (U_i \times U_j)\}_{(i, \,j)
 \in \cc{U}_1}.
$$ 

\vspace{1ex} \item \label{upsilon}
 The morphism $\upsilon$ ``forgets'' the family and the
 trivialization, is given by
$$ 
\upsilon^{*}(X,\, S \to I, \, \theta) \,=\, X,
\hspace{3ex} \upsilon^{*}\; is \; a\;faithful\; functor.
$$

\vspace{1ex} \item \label{varrho}
 The morphism $\varrho$ ``forgets'' the object and the
 trivialization, is given by
$$ 
\varrho^{*}(X,\, S \to I, \, \theta) \,=\, S \to I,
\hspace{3ex} \varrho^{*}\; is \; a\;faithful\; functor.
$$

\vspace{1ex} \item \label{D_U} 
An object of $\cc{D}_\cc{U}$ is a pair $(S \to I, \, \lambda)$, where
$\lambda$ is a $\cc{U}_\bullet$-descent datum on the object $S \to I$ in the topos
$\cc{S}_{/I}$. Such a descent datum is an isomorphism in the topos
$\cc{S}_{/ \cc{U}_1 }$, and it consists of the following data:
$$
 bijections \;\; \lambda_{j\,i}: S_i \to  S_j, \; {(i,\, j) \in
  \cc{U}_1},\;\;such\;that \hspace{15ex}
$$
$$
\hspace{25ex} \lambda_{i\,i} = id, \; \; i \in I, \; \;
 \lambda_{k,i} = \lambda_{k\,j} \circ \lambda_{j\,i}, \;\; 
 (i,\,j,\,k) \in \cc{U}_2.
$$     

\vspace{1ex} \item \label{delta}
 The morphism $\delta $ ``forgets'' the descent datum, is given by
$$ 
\delta^{*}(S \to I, \, \lambda) \,=\, S \to I,
 \hspace{3ex} \delta^{*}\; is \; a\;faithful\; functor.
$$

\vspace{1ex} \item \label{P_U}
An object in the image of the functor
$\rho^*$ is an object in the topos $\cc{E}_{/U}$ of the form
$ \gamma^{*}S \, \times_{\gamma^{*}I} \, U  \mr{} U\, , \;
\{\gamma^{*}S_i \times U_i \mr{} U_i\}_{i
 \in I}$. A $U_\bullet$-descent datum $\sigma$ on such object is an
isomorphism in the topos $\cc{E}_{/  U \times U}$, and it  consists of a
family of isomorphisms in $\cc{E}$,  
$\{\sigma_{j,\,i}\}_{(i, \, j) \in \cc{U}_1}$:
%
%
$$
\xymatrix
         {
            \gamma^{*}S_i \times U_i \times U_j  \ar[r]^{\sigma_{j,\,i}}
                                                 \ar[d]
          & \gamma^{*}S_j \times U_j \times U_i  \ar[d] 
          \\
            U_i \times U_j  \ar[r]^{\tau}  
          & U_j \times U_i
         }
$$
satisfying the apropiate identity and cocycle conditions (the arrow
$\tau$ is the symmetry isomorphism). The commutativity of the square
implies that $\sigma_{j,\,i}$ is completely determined by its first
projection. We abuse notation and write this projection with the same
letter: $\gamma^{*}S_i \times U_i \times U_j  \mr{\sigma_{j,\,i}}
\gamma^{*}S_j$.

\vspace{1ex} \item \label{DtoP}
Given a descent datum as in (\ref{D_U}), the morphism of simplicial
topoi $\rho_\bullet$
  induces a descent datum as in (\ref{P_U}):
$$ 
\sigma = \rho^*_1 \lambda\,, \;\;\;\;
\sigma_{j,\,i} = \gamma^*\lambda_{j,\,i} \times \tau
$$

\vspace{1ex} \item \label{t=dd}
The fact that $\cc{E}_{/U} \mr{\varphi} \cc{E}$ is of effective
descent means (in particular) that given a descent datum as in
(\ref{P_U}), there exist a (unique up to isomorphism) object $X$ in
$\cc{E}$ together with an isomorphism  ${\gamma}^{*}S
\times_{{\gamma}^{*}I} U \mr{\theta}  X \times U$ \mbox{over $U$},
$\{\gamma^{*}S_{i} \times U_{i} \mr{\theta_{i}}  X \times
U_{i}\}_{i \in I}$ 
 (Thus, $X$ is a 
locally constant object $\cc{U}$-split by a trivialization
$\theta$). Moreover, this isomorphism $\theta$ is compatible with the
descent datum $\sigma$ and the trivial descent datum on $X \times
U \to U$, \mbox{$\{X \times U_i \to U_i \}_{i \in I}$,} given by $X \times U_i \times U_j \mr{X \times \tau} X \times U_j
\times U_i$. That is:
$$
\xymatrix
         {
            \gamma^{*}S_i \times U_i \times U_j 
                               \ar[r]^{\sigma_{j,\,i}}
                               \ar[d]^{\theta_i \times U_j}
          & \gamma^{*}S_j \times U_j \times U_i  
                               \ar[d]^{\theta_j \times U_i} 
          \\
            X \times U_i \times U_j  \ar[r]^{X \times \tau}  
          & X \times U_j \times U_i
         }
$$ 

The object $X$ corresponding to a descent datum as in (\ref{DtoP})
furnish the inverse image for the morphisms $\cc{E} \to 
\cc{D}_\cc{U}$ and $ \cc{P}_\cc{U} \to \cc{D}_\cc{U}$.
\end{enumerate}

\vspace{1ex}

\begin{proposition} \label{equivalence}
The push-out topos $\cc{P}_\cc{U}$ of locally constant objects split
by $\cc{U}$ (definition \ref{pushout}) is equivalent to the following
category:

\emph{Objects:} $(S \to I,\, \sigma)$, is a pair, where $S \to I
\, \in \, 
\cc{S}_{/I}$, and $\sigma$ is a $U_\bullet$-descent datum on  
$\gamma^{*}S \, \times_{\gamma^{*}I} \, U  \mr{} U$ in $\cc{E}$
(\ref{facts} (\ref{P_U})).

\emph{Arrows:} $(S \to I,\, \sigma) \to (T \to I,\, \eta)$, is a
family of functions $\{S_i \mr{\vartheta_i} T_i\}_{i \in I}$ compatible with
  the descent data:
$$
\xymatrix
         {
          \gamma^{*}S_i \times U_i \times U_j
	         \ar[r]^{\sigma_{j,\,i}}
                 \ar[d]^{\gamma^{*}\vartheta_i \times U_i \times U_j}
          &
          \gamma^{*}S_j \times U_j \times U_i
                 \ar[d]^{\gamma^{*}\vartheta_j \times U_j \times U_i}
          \\ 
          \gamma^{*}T_i \times U_i \times U_j  
                 \ar[r]^{\eta_{j,\,i}}
          &
          \gamma^{*}S_j \times U_j \times U_i
         }
$$
\end{proposition}
\begin{proof}
We have a functor defined by the assignment \mbox{$(X, S \to I, \theta)
\mapsto (S \to I, \sigma)$,} where the descent data $\sigma_{j,\,i}$
is given by the composite: 
$$
\gamma^{*} S_i
  \times U_i \times U_j 
\;\;\mr{\theta_{i}\times U_j}\;\; X \times U_i \times U_j  \;\;\mr{X 
  \times \tau} \;\;   
X \times U_j \times U_i 
\;\;\mr{\theta_{j}^{-1} \times U_i}\;\;\gamma^{*} S_j \times U_j \times U_i
$$
The statement in \ref{facts} (\ref{t=dd}) says that this functor is an
equivalence of categories. 
\end{proof}

All the details in the proof of this proposition can be checked in an
straightforward way, and it is interesting to do so to understand
exactly how the two types of data  match.

\begin{sinnadaitalica} \label{versus}
When dealing with locally constant objects we shall use the
trivialization or the descent data indistinctly. We shall write $X =
(X, S \to I, \theta)$, or $X = (S \to I, \sigma)$, indicating as usual
with $X$ either the object $X$ or the whole structure.
\end{sinnadaitalica}

 \pagebreak

\begin{facts} \label{facts2}
\end{facts}
\begin{enumerate}
\item \label{surjective}
 The morphism $\cc{S}_{/I} \mr{\varrho} \cc{P}_\cc{U}$ ``forgets'' the descent datum, is given by
$$ 
\varrho^{*}(S \to I, \, \sigma) \,=\, S \to I.
$$

Clearly it is surjective (actually a surjective family of points).
\vspace{1ex} \item \label{connected2}
Given a set $T \in \cc{S}$, the inverse image of the structure
morphism $\cc{P} \mr{\gamma} \cc{S}$ in $T$ is the trivial descent
datum on the constant family $\gamma^*T \times U \to U$, given by $\gamma^*T \times U_i \times U_j \mr{\gamma^*T \times
  \tau} \gamma^*T \times U_j \times U_i$. Thus, $\gamma^*T =
(\gamma^*T,\, \gamma^*T \times  \tau)$. 

\vspace{1ex}

Clearly, if $\cc{E}$ is connected, so it is $\cc{P}$.
\qed 
\end{enumerate}

%
%

\vspace{2ex}


\section{Covering projections associated to a (fixed) cover}
\label{GG_U} 


Let $X = (S \to I, \, \sigma)$ be a locally constant object trivialized by
a cover $U \to \gamma^* I$. 
\begin{definition} \label{actiontriple}
An \emph{action triple} for $X$ is a 3-tuple $(u,\,v,\,s)$ where \mbox{$C
\mr{u} U_i$,} $C \mr{v} U_j$,  $C \in \cc{E}$, $C \neq \emptyset$, and
$S_i \mr{s} S_j$ a bijective function, such that
$$
\xymatrix
         {
            \gamma^{*}S_i \times U_i \times U_j 
                               \ar[r]^{\sigma_{j,\,i}}
          & \gamma^{*}S_j \times U_j \times U_i   
          \\
             \gamma^{*}S_i \times C  
                               \ar[u]_{\gamma^{*}S_i \times (u,\, v)}
                               \ar[r]^{\gamma^* s \times C}  
          &  \gamma^{*}S_j \times C
                               \ar[u]_{\gamma^{*}S_j \times (v,\, u)}
         }
$$
\end{definition}
\begin{remark} \label{remarkunique}
Notice that $C \neq \emptyset$ implies that for given $(u,\,v)$, if
	 there exists a bijection $s$ to complete an action triple,
	 this bijection  is
	 unique.
\qed \end{remark}
\begin{remark} \label{remark2cell}
Given an action triple  $(u,\,v,\,s)$ and an arrow $C' \mr{f} C$, $C'
\neq \emptyset$, 
the pair $(u',\, v')$, where $u' = uf,\, v' = vf$, can be completed
into an action triple $(u',\,v',\,s')$, with $s' = s$. 
\qed\end{remark}
\begin{remark} \label{remarkconnected}
Any pair  $(u,\,v)$ with $C$ a connected object can always be
completed by a bijection $s$ into an action triple.
\qed\end{remark}
The proof of the following proposition is rather straightforward, and
it is left to the reader.
\begin{proposition} \label{epimono}
Let $X \to Y$,  $(S \to I, \sigma) \mr{\vartheta} (T \to
I, \eta)$,  $S \mr{\vartheta} T$ (see \ref{equivalence}) be a morphism between
locally constant objects. Then:

(a) If $\vartheta$ is surjective,  given any action triple $(u,\,v,\, s)$, $S_i \mr{s}S_j$ for $X$, the pair $(u,\,v)$ can be completed into 
an action triple $(u,\,v,\, t)$, $T_i \mr{t} T_j$ for $Y$, $t\vartheta_i = \vartheta_js$.

(b) If $\vartheta$ is injective,  given any action triple $(u,\,v,\, t)$, $T_i
\mr{t} T_j$ for $Y$, the pair $(u,\,v)$ can be completed into  
an action triple $(u,\,v,\, s)$, $S_i \mr{t} S_j$ for $X$, $t\vartheta_i =
\vartheta_js$.

\qed\end{proposition}
\begin{definition} \label{equivalenced} 
The descent data $\sigma$ determines an equivalence relation on the
set $S$ as follows: Given $x, \; y \in S, \; x \in S_i, \; y \in S_j$,
then $x \sim_\sigma y$ if there exists an action triple $(u,\,v,\,s)$
such that $y = sx$. That is,  
$$
\xymatrix@C=15pt
         {
          \gamma^{*}S_i \times C  \ar[rr]^{\gamma^* s \times C}  
          &&  \gamma^{*}S_j \times C
          \\
          & C \ar[ul]^{\gamma^*x  \times  C} \ar[ur]_{\gamma^*y \times C}          
         }
$$
\end{definition}
\begin{remark} \label{remarkmono}
Notice that if $x \sim_\sigma y$, we can always choose an action triple such that the arrow
$C \mr{(u, \, v)} U_i \times U_j$ is a monomorphism.
\qed \end{remark}

This relation is reflexive (given $x \in U_i$, take $C = U_i$, $u =
id$, $v = id$, so that $(u,\,v) = \Delta$, and $s = id$. This
establishes $x \sim_\sigma x$) and clearly symmetric. Its transitive
clousure is an equivalence relation, that we denote also by
$\sim_\sigma$. The \emph{generating pairs} are those
pairs  $x \sim_\sigma y$ which are related by a single action triple.

\begin{proposition} \label{preboolean}
Let $X = (S \to I, \, \sigma)$ be any locally constant object, and $R
\subset S$ be an equivalence class of $\sim_\sigma$.  Given any subobject 
 $Y = (T \to I, \, \sigma)$, $T \subset S$, if
$R \cap T \neq \emptyset$, then $R \subset T$.
\end{proposition}
\begin{proof}
The proof is immediate (consider proposition \ref{epimono} (b)) 
\end{proof}

\begin{corollary} \label{connected}
A locally constant object $X = (S \to I, \sigma)$ is a
connected object in $\cc{P}_\cc{U}$  if and only if $\sim_\sigma $ has
a single equivalent class (that is, $x \sim_\sigma y$ for all pairs $x,\, y$).\qed\end{corollary}
\emph{Warning}: In general the descent datum does not restrict to $R$,
so that equivalent classes are not subobjects.

\vspace{1ex}

The aim of the following development is to prove that the topos of
covering projections (to be introduced in the next section) has generators.

\vspace{1ex}

We shall not define what do we
mean by \emph{pregroupoid}. Here it can be considered to be just a set.


\begin{sinnadastandard} \label{pregroupoid}
{\bf The pregroupoid $\bb{G}_\cc{U}$}
 
A \emph{premorphism} $i \mr{\phi} j$ is  sequences of spans: $i_0 = i,
\; i_n = j,\; n > 0$,  
$$
\phi = ((U_{i_0} \ml{u_0} C_0 \mr{v_0} U_{i_1}), \; (U_{i_1} \ml{u_1} C_1 \mr{v_1}
U_{i_2}), \;\; \ldots \;\; (U_{i_{n-1}} \ml{u_{n-1}} C_{n-1}
\mr{v_{n-1}}  U_{i_n}))
$$
 $(i_k,\,i_{k+1}) \in \cc{U}_1$, and 
$C_k \mmr{(u_k, \,  v_k)} U_{i_k} \times U_{i_{k+1}},
 \; C_k \neq \emptyset$, a subobject of  $U_{i_k} \times U_{i_{k+1}}$.\end{sinnadastandard}
%
%
%
%
We denote the set of premorphisms 
$
\xymatrix
        {
         \bb{G}_\cc{U} \ar@<1ex>[r]^{\delta_0} 
                       \ar@<-1ex>[r]^{\delta_1} 
         & I
        }
$
, $\;\;\delta_0 \phi  = i, \; \delta_1 \phi = j$.

\begin{sinnadastandard} \label{preaction}
{\bf The partial action of  $\bb{G}_\cc{U}$}

Given any locally constant object $X = (S
\to I, \, \sigma)$, 
$
\xymatrix
         {
          \bb{G}_\cc{U} \ar@<1ex>[r]^{\delta_0} 
                        \ar@<-1ex>[r]^{\delta_1} 
         & I
         }
$ 
has a \emph{partial}
\emph{action} $S_i \times \bb{G}_\cc{U}[i, \, j] \mr{s}
S_j$ on the family $(S \to I)$. Given  $i \mr{\phi} j$, the action $S_i
\mr{s_\phi} S_j$ is defined if $\phi$ can be completed
into a sequence of action triples $((u_0,\, v_0,\, s_0),\; (u_1,\,
v_1,\, s_1),\; \ldots \; (u_{n-1},\, v_{n-1},\, s_{n-1}))$. If such
is the case,   
set \mbox{$s_\phi =  s_{n-1}\ldots s_1 s_0$.} 
\end{sinnadastandard}
By definition and remark
\ref{remarkmono} we have:
\begin{proposition}
Given any $x
\in S_i$, $y \in S_j$, 

$\hspace{17ex} x \sim_\sigma y \; \iff \; \exists \, i \mr{\phi}
j \; \in \; \bb{G}_\cc{U} \; | \; s_\phi x = y.$ 
\qed\end{proposition}
We see that a locally constant object is connected if and only if  the action is
transitive (see corollary \ref{connected}). 
Thus all possible families $S \to I$ which admit a connected 
descent datum are quotients of subsets of
$\bb{G}_\cc{U}$. It follows
 
\begin{proposition} \label{generators1}
The collection of all connected locally constant objects trivialized by a
cover $U \to \gamma^*I$ is a (small) set.
\qed\end{proposition}

\vspace{2ex} 
{\bf Covering projections}

Not all locally constant objects should be considered as covering
projections. However, to require that the $U_\bullet$-descent datum
comes from a $\cc{U}_\bullet$-descent datum (\ref{facts} (\ref{DtoP})), as in
the case of topological spaces (section \ref{cpts}), is too
restrictive when the topos is not spacial. 

\vspace{1ex}

\begin{definition} \label{coveringprojectionobject}
We say that a locally constant object $X = (S \to I, \, \sigma)$
trivialized by a cover $U \to \gamma^* I$ is a \emph{covering
  projection} if, for each $(i, \, j) \in \cc{U}_1$, the family
$C \mr{(u, \, v)} U_i \times U_j$ is an epimorphic family, where 
$(u, \, v)$ ranges over all action triples.
\end{definition}
 The following two propositions follow immediately from proposition
 \ref{epimono}.
\begin{proposition}\label{subobjects}
Any subobject in  $\cc{P}_\cc{U}$ of a covering projection is a
covering \mbox{projection.}
\qed\end{proposition} 
\begin{proposition}\label{quotients}
Any quotient in  $\cc{P}_\cc{U}$ of a covering projection is a
covering \mbox{projection.}
\end{proposition}

\begin{proposition}\label{finitelimits}
Any finite limit in $\cc{P}_\cc{U}$ of covering projections is a
covering \mbox{projection.}
\end{proposition}
\begin{proof}
The terminal object clearly is a covering projection. Let
$(u,\,v,\,s)$, $(u',\,v',\,s')$ be action triples for $X = (S 
\to I, \, \sigma)$,  \mbox{$X' = (S' \to I', \, \sigma ')$.} The fiber
product of $C$ and $C'$ over $U_i \times U_j$ (if non 
empty) is an action triple for $X$ and $X'$ simultaneously. By 
construction of binary \mbox{products} in $\cc{P}_\cc{U}$ it readily follows
that it  determines, together
with $s \times s'$, an action triple for the product $X \times
X'$. The proof follows from this and the fact that in a topos the
fiber products of epimorphic families yield an epimorphic family. The
case of a general finite limit can be treated exactly in the
same way, but it follows anyway  from proposition \ref{subobjects}.   
\end{proof}

\begin{proposition} \label{boolean}
Let $X = (S \to I, \, \sigma)$ be a covering projection, and $R \subset S$
be an equivalence class of $\sim_\sigma$. Then, the descent datum
$\sigma$ restricts to $R \to I$ and it determines a subobject 
$A \mmr{} X$, $A = (R \to I, \, \sigma)$. Furthermore, given any subobject 
 $Y = (T \to I, \, \sigma)$, $Y \mmr{} X$, $T \subset S$, if $A \cap Y
\neq
\emptyset$,  (equivalently $R \cap T \neq \emptyset$),  then $A \mmr{}
T$, (equivalently $R \subset T$).
\end{proposition}
\begin{proof}
Let $(i,\,j) \in \cc{U}_1$. Take an epimorphic family $C \mr{(u, \, v)} U_i \times U_j$, with $(u,\, v,\, s)$ action
triples, and  consider the following diagram;
$$
\xymatrix@C=15pt
         {
          \gamma^{*}S_i \times U_i \times U_j 
                               \ar[rr]^{\sigma_{j,\,i}}
          && \gamma^{*}S_j \times U_j \times U_i   
          \\
            \gamma^{*}R_i \times U_i \times U_j
                               \ar@{_{(}->}[u] 
                               \ar@{-->}[rr]^{\sigma_{j,\,i}}
          && \gamma^{*}R_j \times U_j \times U_i
                               \ar@{_{(}->}[u]   
          \\
            \gamma^{*}R_i \times C  
                               \ar[u]_{\gamma^{*}R_i \times (u,\, v)}
                               \ar[rr]^{\gamma^* s \times C}  
          &&  \gamma^{*}R_j \times C
                               \ar[u]_{\gamma^{*}R_j \times (v,\, u)}
         }                              
$$
 The family 
$\xymatrix@1@C=50pt{\gamma^* R_i \times C \;\;  \ar[r]^-{\gamma^* R_i \times (u,
    \, v)}   & \;\;\gamma^{*}R_i 
  \times U_i \times U_j }$ is epimorphic, and the outer diagram
 commutes by definition. It follows that $\sigma_{j,\,i}$ factors as
 shown. The 
 second assertion is immediate (see proposition \ref{preboolean}).
\end{proof}

\begin{proposition} \label{locallyconnected}
Let $(S \to I, \, \sigma)$ be a covering projection. The quotient set $S \to S_{/ \sim_\sigma}$ has the following property:

Given any set $T$,
$T \in \cc{S}$, and a function $S \mr{\vartheta} T \times I$ over $I$ (that is
, a function $S \mr{\vartheta} T$, $\{S_i \mr{\vartheta_i} T\}_{i \in I}$, then $\vartheta$
determines a morphism $(S, \,\sigma) \mr{\vartheta} 
(\gamma^*T,\, \gamma^*T \times  \tau)$ in $\cc{P}_\cc{U}$ if and only
if it factors $S \mr{} S_{/ \sim_\sigma} \mr{} T$. 
\end{proposition}
\begin{proof}
Let $x \in S_i$, $y \in S_j$ be such that $x \sim_\sigma
y$ by an action triple $(u, \, v,\, s)$. Consider the following diagram: 
$$
\xymatrix@C=15pt
         {
          \gamma^{*}T \times U_i \times U_j 
                               \ar[rr]^{\gamma^* T \times \tau}
          && \gamma^{*}T \times U_j \times U_i   
          \\
            \gamma^{*}S_i \times U_i \times U_j
                               \ar[u]_{\gamma^*\vartheta_i \times U_i \times U_j} 
                               \ar[rr]^{\sigma_{j,\,i}}
          && \gamma^{*}S_j \times U_j \times U_i
                               \ar[u]_{\gamma^*\vartheta_j \times U_j \times U_i}   
          \\
            \gamma^{*}S_i \times C  
                               \ar[u]_{\gamma^{*}S_i \times (u,\, v)}
                               \ar[rr]^{\gamma^* s \times C}  
          &&  \gamma^{*}S_j \times C
                               \ar[u]_{\gamma^{*}S_j \times (v,\, u)}
          \\
          & C \ar[ul]^{\gamma^*x \times  C} \ar[ur]_{\gamma^*y \times
         C}                } 
$$
(a) The middle square and the lower triangle commute by definition.

\vspace{1ex} 

\noindent (b) The commutativity of the upper square means that $\vartheta$ is a
morphism in $\cc{P}_\cc{U}$. 

\vspace{1ex}

\noindent (c) If $C
\mr{(u, \, v)} U_i \times U_j$ is a monomorphism, the commutativity of
the outer diagram implies that $\vartheta$ factors through
$S_{/\sim_\sigma}$. 
Notice that this together with $C \neq \emptyset$ implies that we have:  
$$(\gamma^*(\vartheta_i x), \, (v, \, u)) = (\gamma^*(\vartheta_j y), \, (v, \, u))
\iff \vartheta_i x = \vartheta_j y.$$ 

\emph{f morphism $\rimply$ $\vartheta$ factors}: Let $x \sim_\sigma y$. Clearly
it is enough to take a generating pair. Furthermore we can
assume that the arrow $(u, \, v)$ is a monomorphism 
\mbox{(remark \ref{remarkmono}).} Since the outer
diagram commutes, it follows $\vartheta_i x = \vartheta_j y$.

\emph{$\vartheta$ factors $\rimply$ $\vartheta$ morphism}: Take an
epimorphic family $C \mr{(u, \, v)} U_i \times U_j$, with 
$(u,\, v,\,s)$ action 
triples. For all $x \in S_i$ let $y = sx$. The family 
$
\xymatrix@1@C=50pt
             {
               C\;\;  \ar[r]^-{\gamma^* x \times (u,\, v)}   
             & \;\;\gamma^{*}S_i \times U_i \times U_j 
             }
$ 
is epimorphic. Since for all $x$, $(u,\,v)$ the outer diagram
commutes, it follows that the upper square commutes. 
\end{proof}
From proposition \ref{generators1} we have, in particular:
\begin{proposition} \label{generators2}
The collection of all connected covering projections trivialized by a
cover $U \to \gamma^*I$ is a (small) set.
\qed\end{proposition} 
 
\vspace{1ex}

Let $\cc{G}_\cc{U} \subset \cc{P}_\cc{U}$ be the full subcategory
whose objects are sums of covering projections trivialized by the
cover $U \to \gamma^* I$. 

From propositions \ref{quotients},
\ref{finitelimits}, \ref{boolean}, \ref{locallyconnected}
\ref{generators2} and facts \ref{facts2} (\ref{connected2}) it follows

\begin{theorem} \label{atomic} $ $

1. The category  $\cc{G}_\cc{U}$ is an atomic (locally connected and
boolean) topos and the full inclusion is the inverse image of a
geometric morphism \mbox{$\cc{P}_\cc{U} \mr{} \cc{G}_\cc{U}$.} If $\cc{E}$ is
connected, so it is  $\cc{G}_\cc{U}$. 

2. The functor $\cc{G}_\cc{U}  \mr{\varrho^*} \cc{S}_I$,
$\varrho^{*}(S \to I, \, \sigma) \,=\, S \to I$ is the inverse image
of a surjective point  $\cc{S}_{/I} \mr{\varrho} \cc{G}_\cc{U}$. 
\qed\end{theorem}

 We abuse notation and write $\cc{E} \mr{\upsilon} \cc{G}_\cc{U}$ 
and $\cc{S}_{/I} \mr{\varrho} \cc{G}_\cc{U}$
for the composites $\cc{E} \mr{\upsilon} \cc{P}_\cc{U} \mr{}
\cc{G}_\cc{U}$  and $\cc{S}_{/I} \mr{\varrho} \cc{P}_\cc{U} \mr{}
\cc{G}_\cc{U}$ respectively. Notice that the point  $\cc{S}_{/I}
\mr{\varrho} \cc{G}_\cc{U}$ can be thought as a family  
$\cc{S} \mr{\varrho_i} \cc{G}_\cc{U}$, $\varrho_{i}^{*}(S \to I, \,
\sigma) \,=\, S_i$, of (enough) points indexed by
the set $I$. 

\vspace{1ex}

Let 
  $I \to \mbf{G}_\cc{U}$
be the localic groupoid of the points  $\cc{S}_{/I} \mr{\varrho}
\cc{G}_\cc{U}$, with (discrete) set of objects $I$
(this groupoid is explicitly constructed in \cite{D2} section 2), and
let \mbox{$\cc{S}_{/I} \mr{} \beta\mbf{G}_\cc{U}$} be its classifying
topos (\cite{JT}, \cite{M3}, \cite{D2}). There is a comparison morphism
$\beta\mbf{G}_\cc{U} \mr{} 
\cc{G}_\cc{U}$ compatible with the respective points. The next theorem
follows from theorem \ref{atomic} by
\mbox{\cite{JT} VIII, 3. Theorem 1} (also \cite{D1} Theorem
8.4, or, explicitly \cite{D2} Theorem 3.6.4). 
\begin{theorem} \label{U_fundamental}
The comparison morphism  $\beta\mbf{G}_\cc{U} \mr{\cong} \cc{G}_\cc{U}$
is an equivalence which identifies the point $\varrho$ with the
canonical point of $\beta\mbf{G}_\cc{U}$.
\qed\end{theorem}
\begin{observation}
Recall that $\mbf{G}_\cc{U}$ may be chosen to be \emph{etal complete}
\cite{M3}. Actually, the construction in \cite{D2} yields an etale
complete localic groupoid.
\end{observation}

The theory presented here generalizes the known theory for the locally
connected case, bringing  a direct proof for the atomicity of the
topos of locally constant objects split by a covering. From remark
\ref{remarkconnected} it immediately follows that for a locally
connected topos every
locally constant object is a covering projection. Also, in this case, the
points $\cc{S}_{/I} \mr{} \cc{P}_\cc{U} = \cc{G}_\cc{U}$ are essential (therefore
representable). This amounts
to the fact that the fibers of an inverse limit in
$\cc{P}_\cc{U}$ 
are the inverse limit of the fibers in $\cc{S}$ (fact that is not true
if the topos is not locally connected). It follows that the
localic groupoid $\mbf{G}_\cc{U}$ is an ordinary discrete groupoid, and
the representation theorem \ref{U_fundamental} can be easily proved
without recourse to Joyal-Tierney results (see \cite{D1},
\cite{D2}). We have

\begin{proposition}
If the topos $\cc{E}$ is locally connected, every locally constant
object is a covering projection. That is, $\cc{G}_\cc{U} =
\cc{P}_\cc{U}$. Furthermore, the points are representable, and the
groupoid $\mbf{G}_\cc{U}$ is an ordinary groupoid in $\cc{S}$.
\qed\end{proposition}

 Our theorem \ref{U_fundamental} for the non locally
connected case appears to be the first genuine application of
\cite{JT} VIII, 3. Theorem 1 in the galois theory of locally constant
objects. It is worth noticing also that it is the first time a non
prodiscrete localic groupoid appears in this theory, as well as an
atomic topos which is not a Galois topos.

\vspace{2ex}

\section{Geometric morphisms induced by cover refinements}
\label{groupoidsection} 

The covers of a topos form a category $Cov(\cc{E})$, taking as arrows
the family morphism. Given two covers $\cc{U} = (U, I, \zeta)$,
$\cc{V} = (V, J, \xi)$,
an arrow is  a pair $U \mr{h} V$, $I \mr{\alpha} J$, making the
following square commutative:
$$
\xymatrix
        {
          U \ar[r]^{h} \ar[d]^{\zeta} 
        & V \ar[d]^{\xi}
        \\
          \gamma^{*} I \ar[r]^{\gamma^*\alpha} 
        & \gamma^{*} J
        }
$$
We say that $\cc{U}$ \emph{refines}  $\cc{V}$, and call
the morphisms \emph{refinements}. In alternative notation, the arrow  $\cc{U} \mr{(h,\,\alpha)} \cc{V}$ corresponds to 
a family \mbox{$h = \{U_{i} \mr{h_i} V_{\alpha(i)}\}_{i \in I}$.} 

Any two covers have a common refinement, namely, the family 
\mbox{$\{U_i \times V_j\}_{(i, \, j) \in I \times J, \; | \; U_i \times V_j
    \neq \emptyset}$.} 
(remark however  that  $Cov(\cc{E})$ is not a filtered category because given two refinements they cannot be further refined to become a single
one). 

Covers do form a filtered poset $cov(\cc{E})$ directed by the existence of a
refinement, and there is a 
functor \mbox{$Cov(\cc{E}) \to cov(\cc{E})$} which identifies the
different refinements.

\vspace{1ex}

We shall denote $\cc{T}op_\cc{S}$ the 2-category of grothendieck topoi
and geometric morphisms. We have:
\begin{proposition} \label{functoriality}
The construction of the atomic topos $\cc{G}_\cc{U}$ of covering
projections (\ref{atomic}) is functorial :  
$Cov(\cc{E}) \mr{\cc{G}} \cc{T}op_\cc{S}$. Given a
refinement $\cc{U} \mr{(h,\, \alpha)} \cc{V}$ in $Cov(\cc{E})$, we
also denote  $\cc{G}_{\cc{U}} \mr{(h,\, \alpha)} \cc{G}_{\cc{V}}$ the
corresponding geometric morphism. The following diagram commutes: 
$$
\xymatrix
        {
          \cc{E} \ar[r]^{\upsilon_\cc{U}} \ar[rd]^{\!\!\!\!\upsilon_\cc{V}}
        & \cc{G}_{\cc{U}} \ar[d]^{(h,\, \alpha)}
        & \cc{S}_{/I} \ar[l]_{\varrho_\cc{U}} \ar[d]^\alpha
       \\
        & \cc{G}_{\cc{V}} 
        & \cc{S}_{/J} \ar[l]_{\varrho_\cc{V}}
        }
$$
\end{proposition}
\begin{proof}
It follows from the universal property of the push-out that there is a
geometric morphism $\cc{P}_\cc{U} \mr{(h,\, \alpha)} \cc{P}_\cc{V}$,
and that the statement of the theorem holds for the push-out construction
$\cc{P}_\cc{U}$ (\ref{pushout}).
Consider the following diagram of topoi and inverse image functors:
$$
\xymatrix
          {
            \cc{E} 
          & \cc{G}_\cc{U} \ar[l]^{\upsilon^*_\cc{U}} \ar@{^{(}->}[r] 
          & \cc{P}_\cc{U} \ar[r]^{\varrho^*_\cc{U}}  
          & \cc{S}_{/I} 
          \\
          & \cc{G}_\cc{V} \ar[ul]^{\upsilon^*_\cc{V}} \ar@{^{(}->}[r]
          & \cc{P}_\cc{V} \ar[r]^{\varrho^*_\cc{V}} \ar[u]^{(h,\, \alpha)^*} 
          & \cc{S}_{/J} \ar[u]^{\alpha^*}
          }
$$
It is clear that the theorem follows if we prove that the inverse
image functor \mbox{$\cc{P}_\cc{V} \mr{(h,\, \alpha)^*}
  \cc{P}_\cc{U}$} sends covering 
projections to covering projections. We need an explicit description
of this functor.

\begin{remark} \label{aid} 
Let \mbox{$(S \to I,\, \sigma) \in
\cc{P}_\cc{U}$} be the value of the functor $(h,\, \alpha)^*$ on an object
$(T \to J,\, \eta) \in \cc{P}_\cc{V}$,  $(S \to I,\, \sigma) = 
(h,\, \alpha)^* (T \to J,\, \eta)$, where $S = \alpha^* T$ (notice
that $S_i = T_{\alpha(i)}$). By construction, there is a commutative diagram:
$$
\xymatrix
     {
      \gamma^* T_{\alpha(i)} \times V_{\alpha(i)} \times V_{\alpha(j)}
      \ar[r]^{\eta_{\alpha(i),\, \alpha(j)}}
    & \gamma^* T_{\alpha(j)} \times V_{\alpha(j)} \times V_{\alpha(i)}
    \\
      \gamma^* S_{i} \times U_{i} \times U_{j}
      \ar[r]^{\sigma_{i,\,j}}  \ar[u]^{id \times h_i \times h_j}
    & \gamma^* S_{j} \times U_{j} \times U_{i}
      \ar[u]^{id \times h_j \times h_i}
     }
$$
\qed\end{remark} \label{inverseimage} 
{\itshape Continuation of the proof.} For each action triple $(x,\,y)$ for $\eta \,$, consider a pull-back
diagram:
$$
\xymatrix@C=40pt
        {
          C \ar[r]-<20pt,0pt>^{(u,\,v)} \ar[d]
        & U_{i} \times U_{j} \ar[d]^{h_i \times h_j}
        \\
          B \ar[r]-<30pt,0pt>^{(x,\,y)}
        & V_{\alpha(i)} \times V_{\alpha(j)} 
        }
$$
It is easy to check with the aid of the diagram in remark \ref{aid}  that the pair
$(u\,v)$ is an action triple for $\sigma$. Since pulling back an
epimorphic family yields an epimorphic family, this finishes the proof.  
\end{proof}

We can think the functor in the previous proposition as a system
$(\cc{G}_\cc{U})_{\cc{U}  \in Cov(\cc{E})}$ of
topoi and geometric morphisms, or as a system of
categories and inverse image functors, indexed by $Cov(\cc{E})$.
Although $Cov(\cc{E})$ is not filtered, 
we can give a simple description of the colimit of the categories.

\vspace{2ex}

\section{The category and the topos of covering projections}
\label{allprojections} 

The objects of this category are pairs $(X,\, \cc{U})$, where $\cc{U}$
is a cover, and
\mbox{$X = (X, S \to I, \theta) = (S \to I, \sigma)$} is a covering
projection (see \ref{versus}) trivialized by $\cc{U}$.

An arrow $(X,\, \cc{U}) \to (Y,\, \cc{V})$ is map $X \mr{f} Y$ in
$\cc{E}$ such that there exist a common refinement 
$\cc{W} \mr{(h,\, \alpha)} \cc{U}$, $\cc{W} \mr{(l,\, \beta)} \cc{V}$,
and an arrow $(h,\, \alpha)^*(X)  \mr{(f,\, \vartheta)} (l,\,
\beta)^*(Y)$ in $\cc{G}_\cc{W}$, $\upsilon^*(f,\, \vartheta) = f$.
From the fact that any two covers have a common refinement it follows
that given two maps $X \mr{f} Y$ and $Y \mr{g} Z$ in $\cc{E}$, if $f$
and $g$ are arrows in $c\cc{G}(\cc{E})$, then so is the composite
$X \mr{gf} Z$. 

\vspace{1ex}

The hom-sets in the category $c\cc{G}(\cc{E})$
are the filtered colimit (actually a filtered union) of the hom-sets
in the categories $\cc{G}_\cc{W}$, indexed by the 
$\cc{W} \geq \cc{U},\; \cc{W} \geq \cc{V}$ in $cov(\cc{E})$. 

Clearly for each $\cc{U} \in Cov(\cc{E})$ there is a faithful functor
$\cc{G}_\cc{U} \mr{\lambda^*_\cc{U}} c\cc{G}(\cc{E})$, and these
functors form a cone for the system  $(\cc{G}_\cc{U})_{\cc{U}  \in
  Cov(\cc{E})}$  (given $X \in
  \cc{G}_\cc{V}$ and a refinement $\cc{U} \mr{(h,\, \alpha)} \cc{V}$,
  the identity map $X \mr{id} X$ establish an isomorphism 
\mbox{$(X, \, \cc{V}) \cong ((h,\, \alpha)^*(X), \, \cc{U})$} in 
$c\cc{G}(\cc{E})$).  

\begin{proposition} \label{catcolimite} $ $

1. The category $c\cc{G}(\cc{E})$ has
finite limits and the functors $\lambda^*_\cc{U}$ preserve
them.

2. The cone $\cc{G}_\cc{U} \mr{\lambda^*_\cc{U}} c\cc{G}(\cc{E})$ is a
colimit of the system of categories and inverse image functors (indexed by the
category $Cov(\cc{E})$).  
\end{proposition}
\begin{proof}
It needs a straightforward but careful verification which is left to
the interested reader.
\end{proof}  
Notice that  there is a functor  $c\cc{G}(\cc{E}) \mr{\upsilon^*}
\cc{E}$ making the following diagram a commutative diagram of faithful
functors: 
\begin{equation} \label{inverse_cG(E)}
\xymatrix
         {
          \cc{G}_\cc{U} \ar[r]^{\lambda^*_\cc{U}}
                        \ar[rd]^{\!\!\!\!\upsilon^*_\cc{U}}
        & c\cc{G}(\cc{E}) \ar[d]^{\upsilon^*}
        \\
        &  \cc{E}
         }
\end{equation}

Consider in the category $c\cc{G}(\cc{E})$ the grothendieck topology
 generated by all the epimorphic families in $\cc{G}_\cc{U}$ for
 all $\cc{U} \in Cov{\cc{E}}$. Notice that this topology is
 subcanonical. The (small) set of covering projections
 corresponding to the connected objects of
 $\cc{G}_\cc{U}$, $\cc{U}$ running over any (small) \emph{cofinal} set of
 coverings (for example, covering sieves, see \ref{coveringsieve} 3.)
 is a ``topologically generating'' family in the sense of \cite{G2} Expose II 3.0.1 (that is,
 every object in $c\cc{G}(\cc{E})$ is covered by objects in the
 family). It follows that the category of sheaves is legitimate and
 that it is a grothendieck topos (\cite{G2} Expose II 4.11). We shall
 denote this topos $\cc{G}(\cc{E})$. There is
 a full and faithful functor $c\cc{G}(\cc{E}) \to
 \cc{G}(\cc{E})$. Clearly the composite functors 
$\cc{G}_\cc{U} \mr{} c\cc{G}(\cc{E}) \to \cc{G}(\cc{E})$ are the
 inverse image functors of geometric morphisms $\cc{G}(\cc{E}) \to
 \cc{G}_\cc{U}$ which determine a cone for the system 
$(\cc{G}_\cc{U})_{\cc{U}  \in Cov(\cc{E})}$ of topoi and geometric
 morphisms. There is a commutative diagram of surjective geometric
 morphisms:
\begin{equation} \label{geometric_G(E)}
\xymatrix
         {
          \cc{G}_\cc{U}
        & \cc{G}(\cc{E})  \ar[l]_{\lambda_\cc{U}}
        \\
        &  \cc{E} \ar[lu]_{\!\!\!\upsilon_\cc{U}}
                   \ar[u]_{\upsilon}
         }                       
\end{equation}
\begin{theorem} \label{toplimit}
The cone  $\cc{G}(\cc{E}) \mr{\lambda_\cc{U}} \cc{G}_\cc{U}$ is a
limit cone in the 2-category $\cc{T}op_\cc{S}$. That is,
$\cc{G}(\cc{E})$ is the inverse limit of the system of topoi and
geometric morphisms $(\cc{G}_\cc{U})_{\cc{U}  \in Cov(\cc{E})}$.
\end{theorem}
\begin{proof}
It follows immediately from proposition \ref{catcolimite} using
\mbox{\cite{G2} Expose IV 4.9.4} (referred to as \emph{the basic theorem concerning
  classifying topos} in \mbox{ \cite{MM} Chapter VIII, 3).}  
\end{proof}
%
%

\vspace{1ex}


We consider now \emph{covering sieves} as a technical tool in order to
exhibit the topos $\cc{G}(\cc{E})$ as the classifying topos of a
progroupoid.

\vspace{1ex}

Consider any topos $\cc{E}$, and let $C_0 \in \cc{S}$ be a (small) set
of generators. 
\begin{assumption}
We assume that  $C \neq \emptyset \;\; \forall C \in C_0$
\end{assumption}


\begin{sinnadastandard} {\bf Covering sieves}
Any (small) set $I$
of objects of $\cc{E}$ (that is, $I \in \cc{S}$, and each $X \in I$ is
an object $X \in \cc{E}$) determines a family 
$\cc{R}(I) = (\Sigma I,\, I, \,\pi)$, where 
$\Sigma I =  \sum_{X \in I} X$, and 
$\pi = {\textstyle \sum_{X \in I} (X \to 1)}$. Thus, $(\Sigma I)_X =
X$. 

Recall that $I \subset C_0$ is a \emph{sieve} if given any
arrow $C \to D$, with $C \in C_0$, if $D \in I$, 
then $C \in I$.  We say that $I$ is a \emph{covering sieve} if $
\Sigma I \mr{} 1$ is an 
epimorphism, that is, if the family $\cc{R}(I)$ is a cover. Covering
sieves form a (small) poset  $sCov(C_0)$
ordered by inclusion.
A family $\cc{U} = (U,\, I,\, \eta)$ determines a sieve  
\mbox{$s(\cc{U}) = \{C \in C_0 \,|\, \exists i \in I \; and \; C \to U_i\}$.}
\end{sinnadastandard}
\begin{proposition} \label{coveringsieve} $\,$ \\ 
\indent
1. The poset $sCov(C_0)$ is filtered.

2. There is a functor $sCov(C_0) \mr{\cc{R}} Cov(\cc{E})$
   defined by $\cc{R}(I) = (\Sigma(I), I, \pi)$.

3. There is a functor $Cov(\cc{E}) \mr{s}
   sCov(C_0)$, and  the composite $\cc{R}(s\cc{U})$ refines $\cc{U}$. In this
   sense, the functor $\cc{R}$ is cofinal.  

4. Given any two sets of generators, $C_0,\; D_0$, there are
 cofinal morphism of posets \mbox{$sCov(C_0) \mr{} sCov(D_0)$,} 
$sCov(D_0) \mr{} sCov(C_0)$.

5. Given any functor $Cov(\cc{E}) \mr{\cc{F}} \cc{X}$ into a
   category $\cc{X}$, it determines a pro-object
   $(\cc{F}_I)_{I \in sCov(C_0)} $,  
   $\cc{F}_I = \cc{F}_{\cc{R}(I)}$, for each set of generators
   $C_0$, and
   all these pro-objects are isomorphic as pro-objects. 
\end{proposition}
\begin{proof} $ $

1. Given two covering sieves $I$, $J$ $\in sCov(C_0)$, the
   intersection sieve $I \cap J$ is also covering. This follows from
   the fact that given $C \in I$, $D \in J$, the product $C
   \times D$ is covered by objects of $C_0$.


2.  Clearly, given $I, \; J \in sCov(C_0)$, if $I
\subset J$, there is a canonical monomorphic refinement 
$\cc{R}(I) \mono \cc{R}(J)$. 

3. That the sieve $s(\cc{U})$ is covering follows because
   every $U_i$ is covered by objects of $C_0$. The rest of the
   statement is clear.

4. Given any covering sieve $I$ in $sCov(C_0)$, it generates a
   sieve in $sCov(D_0)$, 
   $sI = \{D \in D_0 \;|\; \exists D \to C \;\; with\; C \in I\}$,
   which is covering since the objects in $D_0$
   generate. This defines a morphism of posets. The same
   holds in the other direction. It is clear that $ss(I) \subset
   I$. This finishes the proof.

5. Follows from 1. and 4.
\end{proof}

 Here it is important to remark that although the
  refining sieve of a 
   cover  is canonical, \emph{the refinement itself is not}. There is no
   consistent choice of refinements in such a way that
   the maps  $\cc{R}(s\,\cc{U}) \to \cc{U}$ define a transformation
   natural on $\cc{U}$. 

\vspace{1ex}

The functor $\cc{R}$ is not a real cofinal
 functor (in the sense of \mbox{\cite{G2}, Expose I, 8.)} because $Cov(\cc{E})$ is not a
 filtered category. However, for the particular functor
 $Cov(\cc{E})^{op} \mr{\cc{G}} \cc{C}at$, it follows from the
 construction of the colimit in proposition
 \ref{catcolimite}: 
\begin{proposition} \label{sievecolimit}
Given any set of generators $C_0$ of $\cc{E}$, the canonical functor 
$$Colimit_{I \in sCov(C_0)} \cc{G}_{\cc{R}(I)} \mr{} c\cc{G}(\cc{E})$$
is an equivalence of categories.
\end{proposition}
\begin{proof}
This colimit can be constructed in the same way that the category 
$c\cc{G}(\cc{E})$, and it determines a full subcategory. Then, proposition
\ref{coveringsieve}, 3. suffices to show that the inclusion is
essentially surjective.
\end{proof}

Given any set of generators of a topos $\cc{E}$, consider the system
(a pro-object in the 2-category $\cc{T}op_\cc{S}$) 
$(\cc{G}_{\cc{R}(I)})_{I  \in sCov(C_0)}$, and 
the inverse limit of
topoi and geometric morphisms  $Limit_{I \in sCov(C_0)}
\,\cc{G}_{\cc{R}(I)}$ (known to exists, \cite{G2} Expose VI,
8.2.11). This topos does not depend of the chosen set of
generators (proposition \ref{coveringsieve}, 4.).
\begin{theorem} \label{protopos} $ $

1. The restriction of the cone in theorem \ref{toplimit} : 
 $\cc{G}(\cc{E}) \mr{\lambda_{\cc{R}(I)}} \cc{G}_{\cc{R}(I)}$, is a
 limit cone in the 2-category $\cc{T}op_\cc{S}$. That is,
$\cc{G}(\cc{E})$ is the inverse limit of the protopos 
$(\cc{G}_{\cc{R}(I)})_{I  \in sCov(C_0)}$. 

2. There is a canonical natural transformation (in particular a canonical
   morphism of pro-objects) 
$(\cc{S}_{/I} \mr{\varrho_I} \cc{G}_{\cc{R}(I)})_{I \in sCovS(C_0)}$. The topos
   $\cc{G}(\cc{E})$ is equipped
   with a localic point $Sh(L) \mr{\varrho} \cc{G}(\cc{E})$, where $L$ is the
   inverse limit in the category of localic spaces of the discrete
   spaces determined by the sets $I$.
\end{theorem}
\begin{proof} $ $

1. Follows from proposition \ref{sievecolimit}.

2. Follows from proposition \ref{functoriality} and the fact that
   the assignment of the topos of sheaves is a functor that preserves
   all inverse limits (since it has a left adjoint, the localic
   reflection \cite{JT}).
\end{proof}

\vspace{2ex} 

\section{The fundamental progroupoid of a topos}
\label{fundamentalprogroupoid} 

The statements without proof in this section are justified by the yoga
of the theory of classifying topos of localic groupoids established in \cite{M3}.

A \emph{localic progroupoid} $\mbf{G}$ is a proobject $\mbf{G} =
  (\mbf{G}_\alpha)_{\alpha \in \Gamma}$ in the 2-category of localic
  groupoids. There are 
  natural transformations (in particular, canonical morphisms of
  proobjects) $\mbf{I}_\alpha \to \mbf{G}_\alpha$, 
where $\mbf{I}_\alpha$ are the localic spaces of objects.

\begin{definition}
Given a localic progroupoid $\mbf{G}$ as above, its classifying topos
is the inverse limit topos of the classifying topoi $\beta
\mbf{G}_\alpha$,  $\beta \mbf{G} = Limit_\alpha 
\beta \mbf{G}_\alpha$. It is equipped with a localic point $Sh(\mbf{I})
\to \beta \mbf{G}$, where $\mbf{I}$ is the inverse limit of the
localic \mbox{spaces $\mbf{I}_\alpha$}, $\mbf{I} = Limit_\alpha
\mbf{I}_\alpha$. Notice that $Sh(\mbf{I})$ is also the inverse limit
of the topoi  $Sh(\mbf{I}_\alpha)$,  $Sh(\mbf{I}) = Limit_\alpha Sh(\mbf{I}_\alpha)$.  
\end{definition}
When $\mbf{G}$ is an ordinary (strict) progroup, this is exactly the definition
given in \cite{G2} Expose IV 2.7., where the objects of the topos $\beta
\mbf{G}$ are described explicitly.

\vspace{1ex}

Let $\mbf{I} \to \mbf{gG}$ be the inverse limit of the localic groupoids
$\mbf{G}_\alpha$, \mbox{$\mbf{gG} = Limit_\alpha \mbf{G}_\alpha$.} There is a
comparison functor $\beta \mbf{gG} \to \beta \mbf{G}$ (that is $\beta
Limit_\alpha \mbf{G}_\alpha \to  Limit_\alpha  
\beta \mbf{G}_\alpha$). It is an open problem (plausibly with a
negative answer) to know if this is an equivalence. This is related to
the failure or not of the point  $Sh(\mbf{I}) \to \beta \mbf{G}$ to be
of effective descent. Since $Sh(\mbf{I}) \to \beta \mbf{gG}$ is
always of effective descent, we have:
\begin{proposition}
The comparison morphism  $\beta \mbf{gG} \to \beta \mbf{G}$ is an
equivalence if and only if $Sh(\mbf{I}) \to \beta \mbf{G}$ is of
effective descent
\qed \end{proposition} 
The answer is positive in the classical cases corresponding to the Galois
theory of locally connected topoi. Recall that a morphism of
groupoids $\mbf{G} \to \mbf{H}$
is \emph{composably onto} if it is surjective on commutative
triangles (thus, also on arrows and objects) (see \cite{K} 2.7). 
%
%
Extending SGA4 terminology, we say that a progroupoid is \emph{strict}
if the transition morphisms are composably onto. In \cite{K} 4.18. it
is established that the comparison morphism $\beta \mbf{gG} \to \beta
\mbf{G}$ is an equivalence for any strict progroupoid $\mbf{G}$. In
the particular case of 
strict progroups this was first observed in \cite{T}, and later stated
independently in \cite{M2} (where, furthermore,  the equivalence is
proved for localic progroups whose transition morphisms are open
surjections). 
This equivalence  allows to replace strict
progroups by localic prodiscrete groups in the SGA4 Galois
theory of locally connected topos. 
%
%

\vspace{1ex}

Given any topos $\cc{E}$, consider now the system
$(\cc{G}_\cc{U})_{\cc{U}  \in Cov(\cc{E})}$. By the Morita equivalence
for etal complete 
localic groupoids (\cite{B} 2.6, see also \cite{M3} 7.7) it follows:
\begin{proposition} \label{E_fundamental}
The equivalences 
 $\beta\mbf{G}_\cc{U} \mr{\cong} \cc{G}_\cc{U}$ in theorem
\ref{U_fundamental} determine  a system of localic groupoids  
$(\mbf{G}_\cc{U})_{\cc{U}  \in Cov(\cc{E})}$ and a natural
equivalence \mbox{$(\beta\mbf{G}_\cc{U} \mr{\cong} \cc{G}_\cc{U})_{\cc{U}
  \in Cov(\cc{E})}$} of systems.
\qed\end{proposition}

Given any set of
  generators $C_0$, by restriction this determines a localic
  progroupoid, denoted $\pi_{1}(\cc{E})$, 
$\pi_{1}(\cc{E}) = (\mbf{G}_{\cc{R}(I)})_{I  \in sCov(C_0)}$. This
  progroupoid does not depend on the chosen generators, and it is
  defined to be the \emph{fundamental progroupoid } of the topos
  $\cc{E}$. As before, let $L$ be the
   inverse limit in the category of localic spaces of the discrete
   spaces determined by the sets $I$.
\begin{theorem} \label{fundamentaltopos} $ $

Given any topos $\cc{E}$, the topos $\cc{G}(\cc{E})$ of covering
projections is the classifying
topos of the fundamental localic progroupoid $\pi_{1}(\cc{E})$, by an
equivalence   $\beta \pi_{1}(\cc{E}) \mr{\cong} \cc{G}(\cc{E})$ which
identifies the localic points $Sh(L) \mr{} \cc{G}(\cc{E})$,  
\mbox{$Sh(L) \mr{} \beta \pi_{1}(\cc{E})$.}
\end{theorem}
\begin{proof}
It only remains to indicate that the theorem follows immediately from
the given definitions, proposition \ref{E_fundamental} and theorem \ref{protopos}.
\end{proof}

\vspace{2ex}

\section{The representation of torsors}
\label{groupoid} 

In this last section we prove that the fundamental localic progroupoid
$\pi_{1}(\cc{E})$ represents torsors. Given a group $K \in \cc{S}$ and
  a topos $\cc{E}$, recall that the 
category (groupoid) $\cc{T}ors^K(\cc{E})$ of \mbox{$K$-torsors} (see
below) is equivalent to the category
of geometric 
morphisms from $\cc{E}$ to the classifying topos $\beta K$,
$\cc{T}op_\cc{S}[\cc{E}, \beta K] \cong \cc{T}ors^K(\cc{E})$ (\cite{MM}, Chapter
VIII, Theorem 7). We
shall denote $\cc{G}rpd$, $pro\cc{G}rpd$, the 2-categories of
localic groupoids, localic progroupoids, respectively.
\begin{proposition} \label{repre_1} $ $ 

There is an equivalence of categories $\;$ $pro\cc{G}rpd[\pi_{1}(\cc{E}),\, K]
  \cong \cc{T}ors^K(c\cc{G}(\cc{E}))$.
\end{proposition}
\begin{proof}
$$
\xymatrix@R=3pt
        {
 pro\cc{G}rpd[\pi_{1}(\cc{E}),\, K] 
& \; \cong^1  \\ 
 Colimit_{I  \in sCov(C_0)}  \,\cc{G}rpd[\mbf{G}_{\cc{R}(I)}, \, K]  
& \; \cong^2  \\
 Colimit_{I  \in sCov(C_0)} \,\cc{T}op_\cc{S}[\beta\mbf{G}_{\cc{R}(I)},\, \beta K] 
&  \; \cong^3 \\
Colimit_{I  \in sCov(C_0)} \,\cc{T}op_\cc{S}[\cc{G}_{\cc{R}(I)},\, \beta K]       
&  \; \cong^4 \\ 
Colimit_{I  \in sCov(C_0)} \,\cc{T}ors^K(\cc{G}_{\cc{R}(I)})                 
&  \; \cong^5 \\
 \cc{T}ors^K(c\cc{G}(\cc{E})).
        } 
$$
(1) holds by definition of morphisms of progroupoids, (2) by the
	Morita equivalence for etal complete  
localic groupoids (\cite{B} 2.6, see also \cite{M3} 7.7), (3) by
proposition \ref{E_fundamental}, (4) is clear (see \cite{MM}, Chapter
VIII, Theorem 7), and (5) is an standard
property of filtered colimits.   
\end{proof}

Our next task will be to show that there is an equivalence of
categories  $Tors^K(c\cc{G}(\cc{E})) \cong Tors^K(\cc{E})$.

\vspace{1ex}

Given a topos $\cc{E}$, we shall use letters as
variables to describe 
arrows in $\cc{E}$. We shall denote by a central dot $''\pa\,''$ all
group products and group actions.
Given a set $S \in \cc{S}$, and an element $x \in S$, by the letter
$x$ we shall also indicate the corresponding global section
$1 \mr{x} \gamma^*S$ in the topos.  

Given a group $K$ in $\cc{S}$, 
and given \mbox{$x, \, y \in K$,} we set $x/y = x\pa y^{-1}$ and 
$x\setminus y = x^{-1}\pa y$.  Recall that a
\emph{K-torsor} in a 
topos $\cc{E}$ is an \mbox{object} $T \in \cc{E}$, $T \to 1$ epi, and
an action $\gamma^*K 
\times T \mr{} T$ such that the arrow $\gamma^*K
\times T \mr{\varepsilon} T \times T$ defined by
 $\varepsilon(x, \, u) = (x \pa u, \, u)$ is an isomorphism. There is
an arrow $T \times T 
\mr{} \gamma^*K$, $(u,\,v) \mapsto v/u$ defined by $\varepsilon^{-1}(u,\,v)
= (v/u,\, v)$. Thus, $z\pa u = v  \;\siff\; z = v/u$. It immediately
follows the equation $(x \pa u) \,/\, (y \pa u) \;=\; x/y$. 

Clearly any
torsor $T$ determines in a canonical way a locally constant object
 $T = (T,\,K,\, \varepsilon)$ split by the (singleton family) cover $T \to 1$.
\begin{proposition} $ $ \label{torsorsare}

1. Given any K-torsor in $\cc{E}$, the locally constant object $T =
   (T,\,K,\, \varepsilon)$ is a covering projection (definition
   \ref{coveringprojectionobject}), that is, an object in 
$\cc{G}_{\,T}$.

2. The covering projection  $T = (T,\,K,\, \varepsilon)$ is a K-torsor of
   $\cc{G}_{\,T}$ with the same arrow as action. The group product $K
   \times K \to K$ 
   furnish the function that lifts the action $\gamma^*K \times
   T \mr{} T$ into a morphism of covering projections.
\end{proposition} 
\begin{proof}

1.  The corresponding descent data $\gamma^*K \times T \times T
 \mr{\sigma} \gamma^*K \times T \times T$ is described by   
  $\sigma(z, \, u,\, v) = (v\, / \,(z\pa u)\,,\, v,\, u)$ (see
 \ref{equivalence}). Any pair of elements $x, \, y \in K$
 define an arrow $T \mr{(x,\,y)} T \times T$ in the topos, $(x,\,y)(u)
 = (x\pa u,\,y\pa u)$. Let $K \mr{s} K$ be defined by $s(z) = (y/x)/z $
 (with an inverse given by  $h(z) = z \setminus (y/x)$). It is immediate to check that $(x,\,y,\,s)$ is an action
 triple (\ref{actiontriple}). This proves the statement since the
 family of arrows  $T \mr{(x,\,y)} 
 T \times T$, all  $x, \, y \in K$, is an epimorphic family.

2. $\gamma^*K$ as an object of  $\cc{G}_{\,T}$ has the constant split
   structure, and  $\gamma^*K \times T$ the (cartesian) product split structure 
$\gamma^*(K \times K) \times T  \mr{\varepsilon} (\gamma^* K \times T) \times
   T$, 
   described by \mbox{$\varepsilon(x,\,y,\,u) = (x,\,  y\pa u,\, u)$.} It is
   immediate to check that the group \mbox{product} \mbox{$K \times K \to K$}
   furnish the function that lifts the action $\gamma^*K \times
   T \mr{} T$ into a morphism of covering projections.
\end{proof}
Given two torsors $T, H$, an arrow $T \mr{f} H$ in $\cc{E}$
determines a refinement \mbox{$f = (f,\,id_{\displaystyle 1})$} of the respective covers, so
that $H$ can be viewed as an object \mbox{$f^*H \in \cc{G}_{\,T}$.}
\begin{proposition}\label{torsorsare_2}
An arrow in $T \mr{f} H$ in $\cc{E}$ between torsors is equivariant
(that is, it is a morphism of torsors) if and only if the pair
$(f,\,id_K)$  is a morphism  
$T \mr{} f^*H$ of covering projections (that is, a morphism in 
$\cc{G}_{\,T}$). 
\end{proposition}  
\begin{proof}
It is immediate to check (see remark
\ref{aid}) that the
identity function $K \mr{id} K$ lifts $f$ to a morphism of covering
projections if and only if $f$ respects the actions.
\end{proof}
We consider now a cover $\cc{U}$ and a torsor in the topos
$\cc{G}_\cc{U}$. It consists of a covering projection $T = (T, S \to
I, \theta)$, an action  $\gamma^*K \times T \mr{} T$ which comes
together with an action  $\{K \times S_i \mr{} S_i\}_{i \in I}$ in the
topos $\cc{S}_{/I}$, such that $\theta_i(x\pa s,\, u) = x\pa 
\theta_i(s,\,u)$. The torsor in $\cc{S}_{/I}$ is non canonically
(and it seems choice dependent) isomorphic to the
canonical torsor  $\{K \times K \mr{} K\}_{i \in I}$ (\cite{TT} 8.31,
 \cite{G2} Expose IV 7.2.5). Given a section $I \mr{s} S$, it
 determines an isomorphism $\{K \mr{s_i} S_i, \; s_i(x) = x\pa s_i\}_{i
   \in I}$. It also determines a refinement (that we denote
 $\theta_s$) of covers  
$\theta_s = \{\xymatrix@C=8ex{U_i  \ar[r]^{\theta(s_i ,\, - )} & T}\}_{i \in I}$.
\begin{proposition} \label{torsorsare_3}
Given any torsor in $\cc{G}_\cc{U}$ as above, the pair $(id_T,
\, \{s_i\}_{i \in I})$ establishes an isomorphism 
$\theta_s^*(T,\,K,\, \varepsilon) \mr{\cong}  (T, S \to
I, \theta)$ of torsors in $\cc{G}_\cc{U}$. 
\end {proposition}
\begin{proof}
Notice that the action in both covering projections is the same. It
is immediate to check (see remark
\ref{inverseimage}) that
$\{s_i\}_{i \in I}$ lifts $id_T$ to a morphism of covering projections.
\end{proof}
\begin{theorem} \label{repre_2}
The faithful functor $c\cc{G}(\cc{E}) \mr{\upsilon^*} \cc{E}$
(\ref{inverse_cG(E)}) 
establishes an equivalence of categories $Tors^K(c\cc{G}(\cc{E}))
\mr{\cong} Tors^K(\cc{E})$.
\end{theorem}
\begin{proof}
It follows from propositions \ref{torsorsare_2} and \ref{torsorsare_3}
that the torsor defined in
\mbox{proposition \ref{torsorsare}} determines  an
inverse $Tors^K(\cc{E}) \mr{\cong}  Tors^K(c\cc{G}(\cc{E}))$.
\end{proof}
It follows then from this theorem and proposition \ref{repre_1}
\begin{theorem}\label{representstorsors}
Given any topos $\cc{E}$ and group $K \in \cc{S}$, the fundamental
localic progroupoid 
$\pi_{1}(\cc{E})$ represents $K$-torsors. That is, there is an
equivalence of categories 
$pro\cc{G}rpd[\pi_{1}(\cc{E}),\, K]
  \cong \cc{T}ors^K(\cc{E})$.
\qed\end{theorem}

Of course, this equivalence induces a bijection between the sets of
equivalence classes of objects.

\end{document}